\documentclass[12pt,reqno,a4wide]{amsart}

\usepackage{tikz} 
\usepackage{calrsfs}
\usepackage[charter]{mathdesign}

\allowdisplaybreaks

\title{Elementary evaluation of a determinant of Kirschenhofer and Thuswaldner}

\author[H.~Prodinger]{Helmut Prodinger}
\address{H.~Prodinger\\Department of Mathematical Sciences, Mathematics Division\\ Stellenbosch University, Private Bag X1, 7602 Matieland, South Africa}
\email{hproding@sun.ac.za}

\thanks{}
\date{Written March 24, 2018}

\begin{document}
	\maketitle

In a fresh paper~\cite{KT} that I saw yesterday, Kirschenhofer and Thuswaldner evaluated the determinant
\begin{equation*}
D_s=\det \Bigl(\frac1{(2l)^2-t^2(2i-1)^2}\Bigr)_{1\le i, l\le s}
\end{equation*}
for $t=1$. Consider the matrix $M$ with entries $1/((2l)^2-t^2(2i-1)^2)$ where $s$ might be a positive integer or infinity. In \cite{KT}, the transposed matrix was considered, but that is immaterial when it comes to the determinant. 

The aim of this short note is to provide a completely elementary evaluation of this determinant which relies on the LU-decomposition $LU=M$, which is obtained by guessing. The additional parameter $t$ helps with guessing and makes the result even more general. We found these results:
\begin{gather*}
L_{i,j}=\frac{\prod_{k=1}^j\bigl((2j-1)^2t^2-(2k)^2 \bigr)}{\prod_{k=1}^j\bigl((2i-1)^2t^2-(2k)^2 \bigr)}
\frac{(i+j-2)!}{(i-j)!(2j-2)!},\\
U_{j,l}=\frac{t^{2j-2}(-1)^j16^{j-1}(2j-2)!}{\prod_{k=1}^j\bigl((2k-1)^2t^2-(2l)^2 \bigr)
	\prod_{k=1}^{j-1}\bigl((2j-1)^2t^2-(2k)^2 \bigr)}
\frac{(j+l-1)!}{l(l-j)!}.
\end{gather*}

 Note that
\begin{equation*}
\prod_{k=1}^j\bigl((2i-1)^2t^2-(2k)^2 \bigr)=
(-1)^j4^j\frac{\Gamma(j+1-t(i-\tfrac12))}{\Gamma(1-t(i-\tfrac12))}
\frac{\Gamma(j+1+t(i-\tfrac12))}{\Gamma(1+t(i-\tfrac12))}
\end{equation*}
and
\begin{equation*}
\prod_{k=1}^j\bigl((2k-1)^2t^2-(2l)^2 \bigr)=
4^jt^{2j}\frac{\Gamma(j+\tfrac12+\tfrac{l}{t})}{\Gamma(\tfrac12+\tfrac{l}{t})}
\frac{\Gamma(j+\tfrac12-\tfrac{l}{t})}{\Gamma(\tfrac12-\tfrac{l}{t})};
\end{equation*}
using these formul\ae, $L_{i,j}$ resp. $U_{j,l}$ can be written in terms of Gamma functions.

The proof that indeed $\sum_j L_{i,j}U_{j,l}=M_{i,l}$ is within the reach of computer algebra systems. 
An old version of Maple (without extra packages) provides this summation. 
Consequently the determinant is
\begin{equation*}
D_s=\prod_{j=1}^s U_{j,j}.
\end{equation*}
For $t=1$, this may be simplified:
\begin{align*}
D_s&=\frac1{s!}\prod_{j=1}^s\frac{(-1)^j16^{j-1}(2j-2)!(2j-1)!}{\prod_{k=1}^j(2k-2j-1)(2k+2j-1)
	\prod_{k=1}^{j-1}(2j-2k-1)(2j+2k-1)}\\
&=\frac1{s!}\prod_{j=1}^s\frac{16^{j-1}(2j-1)!^2}{(4j-1)!!(4j-3)!! }
=\frac{4^s}{s!}\prod_{j=1}^s\frac{32^{j-1}(2j-1)!^4 }{(4j-1)!(4j-2)! }\\
&=\frac{4^s\,16^{s(s-1)}}{s!^2}\bigg/\prod_{j=1}^s\binom{4j}{2j}\binom{4j-2}{2j-1}
=\frac{4^s\,16^{s(s-1)}}{s!^2}\bigg/\prod_{j=1}^{2s}\binom{2j}{j}\\
&=\frac{16^{s(s-1)}}{s!^2}\bigg/\prod_{j=0}^{2s-1}\binom{2j+1}{j};
\end{align*}
the last expression was given in \cite{KT}.
 We used the notation $(2n-1)!!=1\cdot 3\cdot 5\cdots (2n-1)$.

\end{document}